\documentclass{amsart}

\usepackage[utf8]{inputenc}
\usepackage[english]{babel}
\usepackage{amssymb}
\usepackage[all]{xy}

\usepackage{xcolor}

\makeindex

\newtheorem{theorem}{Theorem}[section]
\newtheorem{thm*}[theorem]{Theorem*}

\newtheorem{lemma}[theorem]{Lemma}

\newtheorem{proposition}[theorem]{Proposition}
\newtheorem{conjecture}[theorem]{Conjecture}
   
\theoremstyle{definition}
\newtheorem{definition}[theorem]{Definition}
\newtheorem{example}[theorem]{Example}
\newtheorem{remark}[theorem]{Remark}   

\newtheorem*{ack}{Acknowledgments}      

\theoremstyle{remark}

\newcommand{\p}[0]{{\mathbb P}}
\renewcommand{\t}[0]{{\mathbb T}}  

\newcommand{\Ii}[0]{{\mathcal I}}
\newcommand{\Oo}[0]{{\mathcal O}} 
\newcommand{\Ss}[0]{{\mathcal S}} 

\newcommand{\rom}[1]{\ \mbox{#1}\ }

\numberwithin{equation}{section}

\begin{document}

\title[Terracini loci]{On the Terracini locus of projective varieties}

\author[E. Ballico, L. Chiantini]{Edoardo Ballico, Luca Chiantini}
\address[E. Ballico]{Dipartimento di Matematica
\\ Universit\`a di Trento\\ Via Sommarive 14\\ 38123 Povo (TN), Italia}
\email{ edoardo.ballico@unitn.it}
\address[L. Chiantini]{Dipartimento di Ingegneria dell'Informazione e Scienze Matematiche
\\ Universit\`a di Siena\\ Via Roma 56\\ 53100 Siena, Italia}
\email{ luca.chiantini@unisi.it}

\begin{abstract} 
We introduce and study properties of the Terracini locus of projective varieties $X$, which is the locus
of finite sets $S\subset X$ such that $2S$ fails to impose independent conditions to a linear system $L$. Terracini loci
are  relevant in the study of interpolation problems over double points in special position, but they also enter naturally
in the study of special loci contained in secant varieties to projective varieties.\\
We find some criteria which exclude that a set $S$ belongs to the Terracini locus. Furthermore, in the case
where $X$ is a Veronese variety, we bound the dimension of the Terracini locus and we determine examples 
in which the locus has codimension $1$ in the symmetric product of $X$.
\end{abstract}
\keywords{Secant varieties. Singular interpolation}
\subjclass{14N07} 
\maketitle

\section{Introduction}\label{intro}

The paper is devoted to the following problem.
 Consider a projective variety $X\subset\p^N$  of dimension $n$, which we assume to be smooth (even if most of our theory could be arranged 
 to the case $X$ irreducible). Let $p_1,\dots,p_r$ be general points in $X$. Then the tangent spaces $T_{X,p_i}$'s are $r$ linear spaces 
 of dimension $n$, and it is natural to expect that their linear span $\tau_{p_1,\dots,p_r}$ has dimension as large as possible, that is
 $\dim(\tau_{p_1,\dots,p_r}) = \min\{N,r(n+1)-1\}$. If this is not the case, the variety $X$ turns out to be quite special. Due to the celebrated 
 Terracini's Lemma (\cite{Terracini11}), the dimension of $\tau_{p_1,\dots,p_r}$, for $p_1,\dots,p_r$ general, equals the dimension of 
 the $r$-th secant variety $S_r(X)$. Varieties $X$ for which $S_r(X)$ has dimension smaller than the expected value $\min\{N,r(n+1)-1\}$
 are called {\it $r$-defective},  in the terminology of \cite{CCi02a}, and they are classified in many cases (the literature is too vast for even a  rough citation). 

Even when $\tau_{p_1,\dots,p_r}$ has the expected dimension for a general choice of the points, there are possibly special sets $S=\{p_1,\dots,
p_r\}$ for which the span of the tangent spaces drops dimension. These subsets are the object of our study, and their locus, considered as a 
locally closed locus in the symmetric product $X^{(r)}$, will be denoted as {\it the $r$-th Terracini locus of $X$}. Of course when the  points $p_i$'s
are linearly dependent, then $S$ belongs to the Terracini locus. Yet, there are cases  in which linearly independent points have
linearly dependent  tangent spaces, even for large $N$, and our attention is specially devoted to these cases.
 
From the point of view of interpolation theory, the Terracini locus is characterized as the locus of finite sets $S$ whose double
 structures do not impose independent conditions to hyperplanes. The notion naturally generalizes to  line bundles $L$ on a variety $X$, 
 and this is the way in which Terracini loci are defined in Section \ref{defs}.
 Our main examples and applications will concern, though, the case of Veronese varieties, where $X$ is a projective space
 and $L$ is the  line bundle $\Oo_X(d)$, for some degree $d$.
 
From the point of view of secant varieties, the notion of Terracini locus is relevant in the following sense. Several problems on secant varieties
 can be better understood if one starts with the {\it abstract secant variety} $ AS_r(X)$, which is the closure in $\p^N\times X^r$ of the set
 of $(r+1)$-tuples $(u,p_1,\dots,p_r)$ such that the $p_i$'s are linearly independent and $u$ belongs to the span of the $p_i$'s
 (see e.g. Part III of \cite{BC}). It is a standard fact that $S_r(X)$ is the image of $ AS_r(X)$ in the projection $\pi$ to the first factor.
 $ AS_r(X)$ is an object  tamer than $S_r(X)$. For instance, when the points $p_1,\dots ,p_r$ are linearly independent smooth points of $X$, 
 then $AS_r(X)$ is smooth of dimension
 $r(n+1)-1$ at $(u,p_1,\dots,p_r)$, and its tangent space projects to $\tau_{p_1,\dots,p_r}$ in $\pi$. Thus, when $X$ is not $r$-defective,
 the Terracini locus determines points $(u,p_1,\dots,p_r)\in AS_r(X)$ where the differential of $\pi$ drops rank. Avoiding these points is relevant 
 for several, even computational, analyses of secant varieties. When $X$ is a Veronese or a Segre variety, understanding the behavior
 of Terracini loci is important in the determination of the rank and the identifiability of (symmetric and non symmetric) specific tensors. We refer to 
 Section \ref{ident}, and also to \cite{COttVan17b}  or \cite{AngeCVan18} for more details. 

Through the paper, we will give sufficient conditions for a subset $S$ to lie outside the Terracini locus. Then, we will produce examples in which we determine
pairs $(X,L)$, with $L$ very ample, in which the Terracini locus is big, in particular it has components of codimension $1$ in $X^{(r)}$. Our examples will be focused
mainly on the case of Veronese varieties, where $X$ is a projective space $X=\p^n$ and $L=\Oo_{\p^n}(d)$, for some degree $d$.
Then we will produce bounds for the dimension of components of the Terracini locus (see Theorem \ref{cc3}). 

The  examples reveal a strict connection between the theory of Terracini loci of Veronese surfaces and the Brill-Noether theory of special linear series on curves 
with a nodal embedding into the projective plane (see e.g. Example \ref{cod1}). 

A final example of points in the Terracini locus of a Segre variety
shows connections between Terracini loci, the identifiability or the algebraic identifiability of points with respect to $X$,
 and the orbits of finite sets under projective transformations. We aim to understand better the connection
in forthcoming researches. 

Altogether, the examples prove that the study of Terracini loci indicates a perspective, and provides a challenge,
for future studies in algebraic geometry.

\begin{ack} 
Partially supported by MIUR and GSAGA of INdAM (Italy).\\
The authors thank the anonymous referees for many important remarks and for pointing out a gap in the proof of Lemma \ref{cut} and in Example \ref{quin1}. 
\end{ack}

\section{Definitions}\label{defs}

In this section we define the Terracini loci of a projective variety $X$, with respect to a line bundle $L$. 
We work in characteristic $0$, and we will assume that $X$ is integral. Call $n$ the dimension of $X$.
\smallskip

Let $S\subset X$ be a $0$-dimensional scheme. Then we denote with $\ell(S)$ the length of $S$.
If $S$ is a (reduced) collection of points, then $\ell(S)$ is the cardinality $\#S$. We also denote with $\Ii_{S,X}$ the ideal
sheaf of $S$ in $X$. We will omit $X$ when $X=\p^n$.

For any point $p\in X_{\mathrm{reg}}$ we denote with $(2p,X)$ the fat point of $X$ which, in any affine
open subset containing $p$, is defined by $m_p^2$, where $m_p$ is the maximal ideal which defines $p$. The length of $(2p,X)$ is $n+1$.
If $D \subset X$ is a divisor smooth at $p$, then $(2p,D)$ is a proper subscheme of $(2p,X)$, of length
$\ell (2p,D) =2n$, and moreover $(2p,X)\cap D =(2p,D)$ (scheme-theoretic intersection).
We will write simply $2p$, instead of $(2p,X)$, if the scheme $X$ over which we take the fat structure is clear.  
For any   reduced finite set $S\subset X_{\mathrm{reg}}$ denote with $(2S,X)$, or simply $2S$,  the union $\cup_{p\in S}(2p,X)$.

\begin{definition}\label{defter}  Let $L$ be a line bundle on $X$. The \emph{$r$-th Terracini locus} of $(X,L)$ is the subset $\t_r(X,L)$ 
of the symmetric product $X_{\mathrm{reg}}^{(r)}$ defined as
\begin{multline*}\t_r(X,L)=\{ S\subset X_{\mathrm{reg}}:  S\rom{is reduced of cardinality} r \\ \rom{and} h^0(\Ii_{(2S,X)}\otimes L)>h^0(L)-(n+1)r \}.
\end{multline*}
By semicontinuity, the Terracini locus is locally closed in $X_{\mathrm{reg}}^{(r)}$. We call  \emph{closed $r$-th Terracini locus} of $(X,L)$ the closure of $\t_r(X,L)$ 
in the Hilbert scheme of all degree $r$ zero-dimensional subschemes of $X$.
\end{definition}

In order to exclude that $ S\subset X_{\mathrm{reg}}$ belongs to the Terracini locus, it is sufficient to check that $h^0(\Ii_{2S}\otimes L) \le h^0(L)-(n+1)r$.
Since $(n+1)r$ is the length of the structure sheaf of $2S$, it is easy to realize that $ S\subset X_{\mathrm{reg}}$ does not belong to the Terracini
locus if $h^1(\Ii_{(2S,X)}\otimes L)=h^1(L)$. Indeed, more in general, we have

\begin{lemma}\label{indepcon} Let $X$ be a smooth projective variety, $L$ a line bundle on $X$ and $Z\subset X$ a $0$-dimensional scheme. We have
$h^0(\Ii _Z\otimes L)=h^0(L) -\ell (Z)$ if and only if $h^1(\Ii _Z\otimes L) =h^1(L)$.
\end{lemma}
\begin{proof}
Consider the exact sequence of sheaves on $X$:
\begin{equation}\label{eqa1}
0 \to \Ii _Z\otimes L \to L \to L_{|Z}\to 0
\end{equation}
Since $L$ is locally free of rank $1$ and $Z$ is supported on a finite set, the sheaves $L_{|Z}$ and $\Oo _Z$ are isomorphic. Thus
$h^0(L_{|Z}) = h^0(\Oo_Z)=\ell (Z)$. Since $Z$ is a zero-dimensional scheme, we have $h^i(\Oo _Z)=0$ for all $i\ge 1$.
Thus the long cohomology exact sequence of \eqref{eqa1} gives $h^i(\Ii _Z\otimes L)=h^i(L)$ for all $i\ge 2$ (which we do not use), and 
\begin{equation}\label{eqa2} \ell (Z)-h^0(L) +h^0(\Ii _Z\otimes L) = h^1(\Ii _Z\otimes L)-h^1(L).
\end{equation}
\end{proof}

\begin{definition} Let $X$ be an integral projective variety, $L$ a line bundle on $X$ and $Z\subset X$ a $0$-dimensional scheme.
We say that $Z$ \emph{imposes independent conditions on $L$} if $h^0(\Ii_Z\otimes L)=h^0(L)-\ell(Z)$. 

By Lemma \ref{indepcon} $Z$ imposes independent conditions on $L$ if and only if $h^1(\Ii_Z\otimes L)=h^1(L)$. 
Since $h^i(\Oo_Z)=0$ for $i>0$, this is equivalent to say that the restriction map $H^0(L)\to H^0(L_{|Z})$ surjects.
\end{definition}

There are trivial cases of the previous definition that we will not consider in the paper.

For instance, assume that the map $\phi_L:X\to \p^N$ associated with $L$ sends $X$ to a variety of dimension $m<n$. In this case, it is well known
that for any point $p$ the scheme $2p$ imposes at most $m+1$ conditions to sections of $L$. It follows that any $S$ of cardinality $r$ belongs
to the Terracini locus $\t_r(X,L)$.

Also, if $r$ is so big that $h^0(L)-(n+1)r<0$, then every $S$ of cardinality $r$ belongs to the Terracini locus $\t_r(X,L)$, for any $L$.

On the other hand, when $\phi_L(X)$ has dimension $n$, then for $r=1$ the Terracini locus is contained in the locus which maps to the
singular locus of $\phi_L(X)$. Similarly, for any $r$ the
Terracini locus $\t_r(X,L)$ contains any $S$ which intersects the inverse image of the singular locus in the map $\phi_L$. Since we are
mainly interested in the case where $L$ is very ample, we will forget the part of the Terracini locus which arises from
the singularities of $\phi_L(X)$. Notice that if $L$ is very ample, then the Terracini locus $\t_1(X,L)$ is empty. 

Our target will be the case of integral embedded varieties with $L=\Oo_X(1)$, so we will forget the case $r=1$ and also assume that
$n$ is the dimension of $\phi_L(X)$.

\begin{remark}\label{caso}
We will consider mainly the case in which $L$ is the line bundle cut on $X\subset  \p^N$ by the system of  hyperplanes, i.e. $L=\Oo_X(1)$. 

In this situation, for any $p\in X_{\mathrm{reg}}$ the hyperplanes $H$ which correspond to sections of $\Ii_{(2p,X)}(1)$ are precisely the hyperplanes tangent to $X$
at $p$. Thus, we can reinterpret the Terracini locus as the set of all $S\subset X_{\mathrm{reg}}$ of cardinality $r$, such that
the span of the tangent spaces to $X$ at the points of $S$ has dimension $< (n+1)r-1.$ 
\end{remark}

In the rest of the paper we will use often the notion of {\it residue} of a scheme with respect to a hypersurface.
The definition of residue and its main properties can be found e.g. in Remark 2.2 of \cite{BallC13}.\\
 If $Z\subset X$ is any scheme and $T$ is a hypersurface of degree $m$, consider the scheme theoretic intersection $Z\cap T$ and call $Res(Z,T)$ the residue
of $Z$ with respect to $T$.\\
The residue is the scheme defined by  $I_Z:F$, where $I_Z$ is the homogeneous ideal of $Z$ and $F$ is an equation for $T$.\\
For any $d$ the multiplication by $F$ determines the exact sequence of ideal sheaves
\begin{equation}\label{resid} 
0\to \Ii _{Res(Z,T)}(d-m)\to \Ii _{Z}(d) \to \Ii _{Z\cap T,T}(d)\to 0 
\end{equation}
where the rightmost sheaf is the ideal sheaf of $Z\cap T$ as a subscheme of $T$. \\
We recall that when $Z=Z'\cup Z''$, with ${Z'}_{\mathrm{red}}\subset T$ and
the support of $Z''$ does not intersect $T$, then $Z\cap T=Z'\cap T$ and $Res(Z,T)=Z'' \cup Res(Z',T)$.\\
We always have  $\ell(Z\cap T)+\ell(Res(Z,T))=\ell(Z)$. It follows that, if $p$ is a smooth point of $T$
and $Z=(2p,X)$, then $Z\cap T=(2p,T)$ and $Res(Z,T)=\{p\}$.\\
We recall also that if $p\in X$ is a smooth point and $L$ is very ample, then $(2p,X)$ imposes independent conditions on $L$,
i.e. $h^1(\Ii_{(2p,X)}\otimes L)=0$.

\section{An additive result}\label{lems}

In this section $X\subset\p^N$ is a integral variety of dimension $n$ and $L$ is a line bundle on $X$. For any $p\in X_{\mathrm{reg}}$ and any finite set $S\subset X_{\mathrm{reg}}$ set $2p:= (2p,X)$ and $2S:= (2S,X)$.

Let $S$ be a reduced $0$-dimensional subscheme of $X_{\mathrm{reg}}$. For any $p\in S$ we have that $S\setminus \{p\}$ does not contain $p$,
and $\ell(S\setminus \{p\})=\ell(S)-1$ and then $\ell(S\cup 2p) =\ell(S)+n$. We can thus define the following condition.

\begin{definition} Let $X$ be an integral projective variety, $L$ a line bundle on $X$ and $S\subset X_{\mathrm{reg}}$ a reduced $0$-dimensional scheme.
We say that 
\begin{center} $(X,L,S)$ satisfies condition $\dagger$\end{center} 
if the following holds: for every $p\in S$ then $S\cup 2p$ imposes independent conditions on $L$, i.e.
$h^0(\Ii_{S\cup (2p,X)}\otimes L)=h^0(L)-\ell(S)-n$.
\end{definition}

\begin{remark}\label{subs}
By Lemma \ref{indepcon}, it turns out that $(X,L,S)$ satisfies condition $\dagger$ if and only if for all $p\in S$ we have  $h^1(\Ii_{S\cup 2p}\otimes L)=h^1(L)$.
In turn, this is equivalent to say that the restriction map $\rho_p:H^0(L)\to H^0(L_{|S\cup 2p})$ surjects for all $p\in S$. Since $L_{|S\cup 2p}$ is isomorphic
to $\Oo_{S\cup 2p}$, in order to check that the restriction map $\rho_p$ surjects it is sufficient to prove that $\dim(Im(\rho_p))\geq \ell(S\cup 2p)=\ell(S)+n$.
\end{remark}

\begin{remark}\label{sottos} If $Z$ imposes independent conditions to $L$, and $Z'$ is a subscheme of $Z$, then $Z'$ imposes independent conditions to $L$.

Indeed we have $H^0(\Ii_Z\otimes L)\subseteq H^0(\Ii_{Z'}\otimes L)$, so that $H^0(L_{|Z'})$ is a quotient of $H^0(L_{|Z})$. 
Hence when $H^0(L)\to H^0(L_{|Z})$ surjects, then also $H^0(L)\to H^0(L_{|Z'})$ surjects.

In particular, if condition $\dagger$ holds for $(X,L,S)$, then for any subset $S'$ of $S$ condition $\dagger$ also holds for $(X,L,S')$.

Another consequence is the following. Assume that $Z$ imposes independent conditions to $L$. Then $h^0(\Ii_Z\otimes L)=h^0(L)-\ell(Z)$ and also
$h^0(\Ii_{Z'}\otimes L)=h^0(L)-\ell(Z')$. If $\ell(Z')<\ell(Z)$ (in particular this holds when $Z$ is reduced and $Z'\neq Z$), then there exists a section
$v\in H^0(\Ii_{Z'}\otimes L)$ which does not belong to $H^0(\Ii_Z\otimes L)$. 
\end{remark}

\begin{lemma}\label{baseloc}
Let $X$ be an integral projective variety, $S\subset X_{\mathrm{reg}}$ a finite set and let $L_1$, $L_2$ be line bundles on $X$. Assume that
 the base locus of $L_2$ does not intersect $S$. If $(X,L_1,S)$ satisfies $\dagger$, then also $(X,L_1\otimes L_2,S)$ satisfies $\dagger$.
\end{lemma}
\begin{proof}
Fix $p\in S$. Since  the base locus of $L_2$ does not intersect $S$ and $S$ is a finite set, then 
there is a section $v\in H^0(L_2)$ which is nowhere zero in an open neighborhood of $S$. The multiplication by $v$ shows
that the rank of the restriction map $H^0(L_1\otimes L_2)\to H^0((L_1\otimes L_2)_{|S\cup 2p})$ cannot be smaller than the rank 
of the restriction map $H^0({L_1})\to H^0({L_1}_{|S\cup 2p})$. The conclusion follows from Remark \ref{subs}.
\end{proof}

\begin{theorem}\label{i1}
Let $X$ be an integral projective variety. Let $L_1, L_2$ be line bundles on $X$ and $S\subset X_{\mathrm{reg}}$ a finite set of length $r$.
If $(X,L_1,S)$ satisfies $\dagger$ and $S$ imposes independent conditions to $L_2$, then $2S$ imposes independent conditions to $L_1\otimes L_2$,
i.e. $ h^0(\Ii_{2S}\otimes L_1\otimes L_2)=h^0(L_1\otimes L_2)-(n+1)r$, so that $S$ is not in the $r$-th Terracini locus of $L_1\otimes L_2$.
\end{theorem}
\begin{proof}
Set $M:= L_1\otimes L_2$. First assume $r=1$ and write $S =\{p\}$. In this case $S\cup 2p=2p$ and $p$ imposes $1$
condition to $L_2$, so that $p$ does not lie in the base locus of $L_2$. Thus $(X,L_1,S)$ satisfies $\dagger$ by assumption, and then
$(X,M,S)$ satisfies $\dagger$ by Lemma \ref{baseloc}.

Now assume $r > 1$. By Remark \ref{sottos} and induction on the integer $r$ we may assume that any subset $S'\subsetneq S$ imposes independent conditions
to $M$. Fix $p\in S$ and take $S'= S\setminus \{p\}$.
We have $h^0(\Ii _{2S'}\otimes M) =h^0(M)-(n+1)(r-1)$. Thus it is sufficient to prove that $h^0(\Ii _{2S}\otimes M) \le  h^0(\Ii _{2S'}\otimes M) -n-1$.

Since $S$ imposes independent conditions to $L_2$ and it is reduced, then by the last part of Remark \ref{sottos} there exists $v\in H^0 (\Ii _{S'}\otimes L_2)$
which does not belong to $H^0(\Ii _S\otimes L_2)$. In particular, $v$ cannot vanish at $p$.
Since by assumptions and by Remark \ref{sottos} $(X,L_1,S')$ satisfies $\dagger$, then $h^0(\Ii _{S'}\otimes L_1))-h^0(\Ii _{S'\cup 2p}\otimes L_1)=n+1$.
Fix  a linear subspace $V\subseteq H^0(\Ii _{S'}\otimes L_1)$ of dimension $n+1$ such that 
$H^0(\Ii_{S'\cup 2p}\otimes L_1)\oplus V = H^0(\Ii _{S'}\otimes L_1)$. Notice that no non-zero elements of $V$ can belong to 
$H^0(\Ii_{2p}\otimes L_1)$, for $H^0(\Ii _{S'}\otimes L_1)\cap H^0(\Ii_{2p}\otimes L_1)= H^0(\Ii_{S'\cup 2p}\otimes L_1)$.
Fix generators $u_0,\dots,u_n$ for $V$. 
Since $v\ne 0$, the sections $u_0\otimes v, \dots u_n\otimes v$ of $H^0(M)$ are linearly independent. Since both $u_i$ and $v$ vanish at each point of $S'$,
we have $u_i\otimes v\in H^0(\Ii _{2S'}\otimes M)$ for all $i$. Hence $V\otimes v$ is a subspace of dimension $n$ in $H^0(\Ii _{2S'}\otimes M)$.
We claim  that $(V\otimes v)\cap H^0(\Ii_{2S}\otimes M)=(0)$. Indeed any non-zero element of $V\otimes v$ is of the form $w\otimes v$, with $w$ non-zero
in $V$. Since $w$ does not belong to $H^0(\Ii_{2p}\otimes L_1)$ and $v$ does not vanish at $p$, then $w\otimes v\notin H^0(\Ii_{2p}\otimes M)$,
hence $w\otimes v\notin H^0(\Ii_{2S}\otimes M)$.
It follows that  $h^0(\Ii _{2S}\otimes M) \le  h^0(\Ii _{2S'}\otimes M) -n-1$.
\end{proof}

For the applications of the following sections, we will need a technical result.

\begin{lemma}\label{cut} 
Let $X$ be an integral projective variety, $S\subset X_{\mathrm{reg}}$ a finite set of cardinality at least $2$ and let $L_1$, $L_2$ be line bundles on $X$. Assume that
for all $p\in S$ there exists a  divisor in $|L_2|$ which is smooth at $p$ and misses the points of $S\setminus p$. 
If $L_1\otimes L_2$ is very ample and $h^1(\Ii_S\otimes L_1)=0$ then $(X,L_1\otimes L_2,S)$ satisfies $\dagger$.
\end{lemma}
\begin{proof}
For any $p\in S$ fix $T\in |L_2|$ which is smooth at $p$ and misses the remaining points of $S$. Then $(S\cup 2p)\cap T=(2p,T)$ and the residue of $S\cup 2p$ 
with respect to $T$ is $S$. Moreover $p$ does not belong to the base locus of $L_1\otimes L_2$, hence it does not belong to the base
locus of $L_1\otimes L_2$ restricted to $T$. Since $L_1\otimes L_2$ is very ample, then $h^1(T,\Ii_{(2p, T)}\otimes (L_1\otimes L_2)_{|T})=0$.
The residue sequence \eqref{resid} becomes
$$0\to \Ii _{S}\otimes L_1\to \Ii _{S\cup 2p}\otimes L_1\otimes L_2 \to \Ii _{(2p,T)}\otimes (L_1\otimes L_2)_{|T}\to 0$$
and gives $h^1(\Ii_{S\cup 2p}\otimes L_1\otimes L_2)=0$.
\end{proof}

\section{Veronese varieties}

In order to understand the range of application of Theorem \ref{i1}, let us illustrate how it applies to some Veronese varieties.

In this section $X$ is the projective space $\p^n$ and we consider the line bundle $L=\Oo_{\p^n}(d)$ and the complete linear system $|L| = |\Oo_{\p^n}(d)|$, so that sections
of $L$ correspond to hypersurfaces of degree $d$. This amounts to be essentially the same as
considering the Veronese variety, image of the Veronese map $v_{n,d}:\p^n\to \p^N$, where $N=\binom {n+d}d-1$, with the linear
system cut by the hyperplanes of $\p^N$. In this section we write $2p$ and $2S$ instead of $(2p,\p^n)$ and $(2S,\p^n)$.

We will mainly focus on the case $n=2$, for small values of $d$.

\begin{example} Assume $n=2, d=3$. Since for $r>3$ one has $h^0(\Oo_{\p^2}(3))-3r<0$, we will consider the cases $r=2,3$.
For $r=2$, the Terracini locus is empty. Namely if $S$ has cardinality $2$, every cubic containing $2S$ also contains the line
joining the two points. Thus $H^0(\Ii_{2S}(3))$ is isomorphic to the space of conics passing through $S$, which has dimension $4=h^0(\Oo_{\p^2}(3))-6$.

For  $r=3$, if $S$ is not aligned, then a cubic containing $2S$ contains the three lines through any pair of points of $S$, so there is just one such cubic and $S$ is
not in the Terracini locus. On the other hand if $S$ is aligned then all cubics through $2S$ contain twice the line defined by $S$, 
so that $h^0(\Ii_{2S}(3))=h^0(\Oo_{\p^2}(1))=3> h^0(\Oo_{\p^2}(3))-9$, thus $S$ belongs to the Terracini locus.
\end{example}

\begin{remark} The Alexander-Hirschowitz theorem \cite{AlexHir95} (we will refer to this theorem with AHThm, in the sequel) guarantees that
there are few cases  of Veronese varieties whose Terracini locus is dense in $X^{(r)}$. They are listed below:
\begin{itemize}
\item $d=2$, any $n>1$, any $r>1$.
\item $d=3$, $n=4$, $r=7$.
\item $d=4$, $n=2$, $r=5$.
\item $d=4$, $n=3$, $r=9$.
\item $d=4$, $n=4$, $r=14$.
\end{itemize}
\end{remark}

\begin{example} Assume $n=2, d=4$. The values of $r$ to consider are $r=2,3, 4,5$.

For $r=2$, since all the subsets of cardinality $2$ in $\p^2$ are projectively equivalent, by the AHThm  the Terracini locus is empty. 

For $r=3$, one shows, as in the previous case, that any set of $3$ aligned points is in the Terracini locus. The Terracini locus contains just the aligned subsets, 
by the AHThm, because all the non-aligned subsets of $3$ points in $\p^2$ are projectively equivalent.

Assume $r=4$. If three points $p_1,p_2,p_3$ of $S=\{p_1,p_2,p_3,p\}$ are contained in a line $Y$, then $S$ belong to 
$\t_4(\p^2,\Oo_{\p^2}(4))$. Namely $h^0(\Oo_{\p^2}(4))-12=3$, while $|\Ii_{2S}(4)|$ contains quartics which split in $Y$ plus a cubic passing 
through $S\cup 2p$. Since $\ell(S\cup 2p)=6$, then  $h^0(\Ii_{S\cup 2p}(3))\geq 4$. All the subsets of $4$ points in $\p^2$, no three aligned, are projectively
equivalent and form a dense subset of $(\p^2)^{(4)}$, thus by the AHThm none of them lies in the Terracini locus.

The case $r=5$ falls in one particular case of the AHThm, so that the Terracini locus is dense in $(\p^2)^{(5)}$.
\end{example}

Some of the previous arguments generalize.

\begin{remark} The case $r=2$ of the two previous examples generalizes. For any $n$ all the subsets of cardinality $2$ in $\p^n$ are projectively equivalent.
Hence, by the AHThm,  the Terracini locus is empty. 
\end{remark}

\begin{proposition}\label{4.5} Let $S$ be a finite reduced subset of $\p^n$ and fix a degree $d$. Assume that $u$ points of $S$ 
belong to a hypersurface $F$ of degree $d-q$, where $u,q$ are positive integers which satisfy 
\begin{equation}\label{trivia}  
\binom {q+n}n +nu > \binom {d+n}n.
\end{equation}
Then $S$ belongs to the Terracini locus $\t_r(\p^n,\Oo_{\p^n}(d))$.
\end{proposition}
\begin{proof} Just an easy count of parameters. Counting hypersurfaces which split in the union of $F$ and a hypersurface of degree $q$
which contains $2(S\setminus F)$ and passes through $F\cap S$, we get:
$$ h^0(\Ii_{2S}(d))\geq  h^0(\Oo_{\p^n}(q)) - (n+1)(r-u)-u, $$
and the last expression is bigger than $h^0(\Oo_{\p^n}(d)) -(n+1)r$, when \eqref{trivia} holds.
\end{proof}

\begin{example}\label{quin1} Assume $n=2, d=5$. This is the first case in which we can appreciate the sharpness of Theorem \ref{i1}.
The values of $r$ that we must consider are $r=3, 4,5,6,7$.

For $r=3$  the Terracini locus is empty. Notice indeed that three double points impose independent conditions to quintics, even when they are aligned.

Assume $r=4$. If  $S$ is aligned, then it belongs to the Terracini locus, because of  Proposition \ref{4.5}. Otherwise, $S$ imposes independent conditions to conics,
i.e. $h^1(\Ii_S(2))=0$. Thus, the assumptions of Lemma \ref{cut} are satisfied  for $L_1=\Oo_{\p^2}(2)$ and $L_2=\Oo_{\p^2}(1)$. This implies that $(\p^2,S,\Oo_{\p^2}(3))$ 
satisfies condition $\dagger$. Hence, by Theorem \ref{i1} $S$ does not belong to the Terracini locus.

Assume $r=5$. If $4$ points of  $S$ are aligned, then $S$ belongs to the Terracini locus, because of  Proposition \ref{4.5}. 
Otherwise, $S$ imposes independent conditions to conics,
i.e. $h^1(\Ii_S(2))=0$. Thus, the assumptions of Lemma \ref{cut} are satisfied  for $L_1=\Oo_{\p^2}(2)$ and $L_2=\Oo_{\p^2}(1)$. This implies that $(\p^2,S,\Oo_{\p^2}(3))$ 
satisfies condition $\dagger$. Hence, by Theorem \ref{i1} $S$ does not belong to the Terracini locus.

Assume $r=6$. If $S$ lies in a conic, then it belongs to the Terracini locus, because of   Proposition \ref{4.5}. Otherwise, $S$ imposes independent conditions to conics,
i.e. $h^1(\Ii_S(2))=0$. Thus, the assumptions of Lemma \ref{cut} are satisfied  for $L_1=\Oo_{\p^2}(2)$ and $L_2=\Oo_{\p^2}(1)$. This implies that $(\p^2,S,\Oo_{\p^2}(3))$ 
satisfies condition $\dagger$. Hence, by Theorem \ref{i1} $S$ does not belong to the Terracini locus.
\end{example}

We can extend some cases of the previous example in the following

\begin{proposition} \label{meta} Assume $S$ is a reduced finite subset of $\p^n$, of cardinality $r$, which satisfies $h^1(\Ii_S(q))=0$ for some $q<d/2$. 
Then $S$ does not belong to the Terracini locus $\t_r(\p^n,\Oo_{\p^n}(d))$.
\end{proposition}
\begin{proof} Take $L_1=\Oo_{\p^n}(q)$ and $L_2=\Oo_{\p^n}(1)$. Then the conditions of Lemma \ref{cut} are satisfied, which implies that $(X,S,\Oo_{\p^n}(q+1))$
satisfies condition $\dagger$. Since $d-q\geq q+1$, then $(X,S,\Oo_{\p^n}(d-q))$ satisfies condition $\dagger$ by Lemma \ref{baseloc}. The conclusion follows
from Theorem \ref{i1}.
\end{proof}

\begin{example} Let us continue with the case  $n=2, d=5$. It remains to analyze what happens for $r=7$.

If $S$  contains  $6$ points on a conic, then it belongs to the Terracini locus, by  Proposition \ref{4.5}. 
This case includes, as a subcase, the situation in which $S$ has $4$ aligned points.

We want to prove that if no $6$ points of $S$ lie in a conic, then $S$ sits outside the Terracini locus. Notice that we cannot hope to do that
by means of Proposition \ref{meta}, because $7$ points cannot impose independent conditions to conics.

Since $r=7$, we have $h^0(\Oo_{\p^2}(5))-3r=0$. Thus we get the claim if we prove that there are no quintics containing $2S$, unless 
$S$ has $6$ points on a conic.

So, assume that $S$ has no $6$ points on a (possibly reducible) conic, and that there exists a quintic $C$ containing $2S$.
Quintics  singular at $7$ points have negative arithmetic genus, thus they are reducible. Assume that $C$ splits in an irreducible
conic plus a cubic. Then two points of $S$ lie outside the conic, thus the cubic must be reducible. It follows that $C$ has one line $Y$ as a component.

Assume that $C$ contains $2Y$. Since $Y$ contains at most three points of $S$, the residual cubic is singular at three points.
No reduced cubic is singular at three points, thus $C$ contains a second double line $2Y'$. As the conic $Y\cup Y'$ cannot contain more than
$5$ points of $S$, the residue of $C$ with respect to $2Y\cup 2Y'$, which is a line, must be singular a two points, a clear contradiction. 
Thus no lines are contained in $C$ twice.

The residue $C'$ of $C$ with respect to $Y$ is a quartic with $4$ singular points passing through the points of $S$, 
hence $C'$ is reducible. If $C'$ is a double conic $2W$, then $W$ contains
all the points of $S$, a contradiction. If $C'$ is the union of two 
irreducible conics $W_1\cup W_2$, then necessarily the $4$ points of $W_1\cap W_2$ belong to $S$ while $Y$ contains the three remaining points of $S$.
Since these three points must be contained in $W_1\cup W_2$, then there exists one of the two conics passing through a total of
$6$ points of $S$, a contradiction. Thus $C'$ contains a line $Y'$. 

$C'$ cannot contain $2Y'$, for what we said above, moreover the conic $Y\cup Y'$ 
contains at most $5$ points of $S$. Thus the residual cubic splits either in one line $Y''$ plus an irreducible conic $C''$, 
or in three lines, since it cannot contain a double line.
In the former case, $Y''\cup C''$ is singular at two points of $S$ and contains at least $4$ points of $S\cap (Y\cup Y')$. If $Y''\cup C''$
contains $5$ points of $S\cap (Y\cup Y')$, then either $C''$ contains $6$ points of $S$ or $Y''$ contains $4$ points of $S$, a contradiction.
If $Y''\cup C''$ contains $4$ points of $S\cap (Y\cup Y')$, then $Y\cap Y'\in S$, and either $C''$ contains $6$ points of $S$ or $Y''$ together
with one of the two lines $Y,Y'$ contain $6$ points of $S$. In both cases we get a contradiction.

Finally, if $C$ splits in $5$ lines, then each point of $S$ belongs to at least two lines of $C$: this is combinatorically  impossible if any line contains at most three  
points of $S$.
\end{example}

In particular, when $L=\Oo_{\p^2}(5)$ on $\p^2$ we see that the Terracini locus is empty for $r=3$, it has codimension $2$ when $r=4,5$,
and codimension $1$ for $r=6,7$.

\section{Veronese varieties and high rank}\label{verhigh}

Recall that, by definition, the Terracini locus is trivial as soon as $h^0(L)-(n+1)r<0$, so the study of the Terracini locus makes sense only for 
\begin{equation}\label{rmax}   r\leq \frac{h^0(L)}{n+1}.
\end{equation}
In this section we  analyze some cases in which $r$ is close to  the upper bound.

\begin{example}\label{sextiche} Consider the Veronese surface with $d=6$, i.e. the case where $L$ is the system of plane sextics. Here the study of the Terracini locus 
makes sense only for $r\leq 9$. Let us examine closely what happens for $r=9$.

Since $h^0(\Oo_{\p^2}(6))-3\cdot 9 =1$, by the AHThm for a general choice of a set of $9$ points $S$ there exists a unique sextic curve with $9$ 
nodes at the points of $S$. In this particular case (which is weakly defective, in the sense of \cite{CCi02a}) the unique sextic with nodes at a 
general $S$ is the double cubic $C$, where $C$ is the unique cubic curve through $S$. 

On the other hand, there are special sets $S$ for which one finds a pencil of sextics passing through $2S$. One case is given by
complete intersections of two cubics.
Indeed if $S$ is contained in two cubic curves, then $h^0(\Ii_{2S}(6))\geq 3$, because $2S$ sits in all pairs of cubics through $S$. 
Notice that complete intersections of this type fill a locus of codimension $2$ in $(\p^2)^{(9)}$, because $8$ points of $X$ can be chosen freely, and 
after that the ninth point is fixed.

We want to prove the existence of other components of the Terracini locus $\t_9(\p^2,\Oo_{\p^2}(6))$ of low codimension in $(\p^2)^{(9)}$. To this aim,
notice that an irreducible sextic curve with $9$ nodes has geometric genus $1$. If $D$ is such a curve and $S$ is the set of nodes of $C$, then 
$2S$ belongs to the base locus of the pencil of sextics generated by $D$ and by $2C$, where $C$ is a cubic curve through $S$. Thus $S$ sits in the Terracini locus.

To compute the dimension of the locus of sets $S$ which are the nodes of irreducible sextics, we proceed as follows.
First observe that $h^0(\Ii_{2S}(6))$ cannot be greater that $2$. Otherwise, fix a pencil of sextics in $H^0(\Ii_{2S}(6))$, disjoint from $D$.
There exists an element $D'$ of the pencil which passes through a general point $P\in D$. Then the length of the intersection between
$D$ and $D'$ is at least $4\ell(S)+1=37$, contradicting Bezout's formula. Then, consider the $1$-dimensional moduli space of elliptic curves. The choice 
of $6$ general points on a general elliptic curve $\Gamma$ determines a linear series of dimension $5$ on $\Gamma$, which contains a family of dimension
$3(6-3)=9$ linear series $g^2_6$. Thus, the pairs $(\Gamma, g^2_6)$ form an irreducible family of dimension $10$. A general $g^2_6$
maps the corresponding curve  $\Gamma$ to an irreducible sextic curve with $9$ nodes $S$. Composing with the automorphisms of $\p^2$
(a $8$-dimensional family) and considering that a general $S$ sits in a pencil of sextics nodal at $S$, we find a $17-$dimensional
family of sets in the Terracini locus $\t_9(\p^2,\Oo_{\p^2}(6))$. Thus the Terracini locus has a component of codimension $1$.
\end{example}

In the following examples we will often use the main result of \cite{Treger89}, which says that, under some numerical assumptions, two
nodal plane curves of the same degree cannot share the same set of nodes.

\begin{example}\label{settiche} Consider the Veronese surface with $d=7$. Here the maximal value $r$ for which the Terracini locus is non trivial if $r=12$.
By the AHThm, a general set $S$ of cardinality $12$ in $\p^2$ is not the set of nodes of septic curves. Let us show the existence of loci
of codimension $1$ in $(\p^2)^{(12)}$ formed by sets of nodes of irreducible septic curves.

Septic curves with $12$ nodes and no other singularities have geometric genus $3$. The moduli space of curves of genus $3$ is $6$-dimensional.
On a general curve of genus $3$ the choice of $7$ general points produces a linear series $g^4_7$, which contains a family of $g^2_7$ of
dimension $3(5-3)=6$. Thus, on a general curve $\Gamma$ of genus $3$ there is a family of series $g^2_7$ of dimension $7-4+6=9$. The general
series $g^2_7$ maps  $\Gamma$ to a plane septic with a set $S$ of $12$ nodes. Such $S$ belongs to the Terracini locus $\t_{12}(\p^2,\Oo_{\p^2}(7))$.
The map from $\Gamma $ to $\p^2$ can be composed with an automorphism of $\p^2$, so we get a family of dimension $9+8=17$ of such maps. 
Thus, the family of pairs (curve $\Gamma$ of genus $3$, birational map $\Gamma\to\p^2$ of degree $7$) is $23$-dimensional and determines a 
family of sets $S$ in the Terracini locus $\t_{12}(\p^2,\Oo_{\p^2}(7))$. By \cite{Treger89}, a general septic of genus $3$ cannot share
the set of nodes with another septic. Thus we have a $23$-dimensional component of the Terracini locus in $(\p^2)^{(12)}$.
\end{example}

\begin{example}\label{ottiche} The previous example can be repeated almost verbatim for $d=8$.
In this case the maximal value $r$ for which the Terracini locus is non trivial if $r=15$.
By the AHThm, a general set $S$ of cardinality $15$ in $\p^2$ is not the set of nodes of octic curves. 

Irreducible octic curves with $15$ nodes and no other singularities have geometric genus $6$. The moduli space of curves of genus $6$ is $15$-dimensional.
On a general curve of genus $6$ the choice of $8$ general points produces a linear series $g^2_8$. Thus on a general curve $\Gamma$ of genus $3$ there
 is a family of series $g^2_8$ of dimension $8-2=6$. The general
series $g^2_8$ maps  $\Gamma$ to a plane octic with a set $S$ of $15$ nodes. Such $S$ belongs to the Terracini locus $\t_{15}(\p^2,\Oo_{\p^2}(8))$.
The map from $\Gamma $ to $\p^2$ can be composed by an automorphism of $\p^2$, so we get a family of dimension $6+8=14$ of such maps. 
Thus, the family of pairs (curve $\Gamma$ of genus $6$, birational map $\Gamma\to\p^2$ of degree $8$) is $29$-dimensional and determines a 
family of sets $S$ in the Terracini locus $\t_{15}(\p^2,\Oo_{\p^2}(8))$. By \cite{Treger89}, a general octic of genus $3$ cannot share
the set of nodes with another octic. Thus we have a $29$-dimensional component of the Terracini locus in $(\p^2)^{(15)}$.
\end{example}

We must stop here the series of examples of codimension $1$ components of the Terracini locus. For a general choice of the degree we cannot guarantee
the existence of such a component.

\begin{example}\label{cod1} For $n=2$ fix a degree $d$ not divisible by $3$ and greater than $8$. Then the maximal $r$ for which the Terracini locus is non trivial is
$$r=\frac{h^0(\Oo_{\p^2}(d))}3 = \frac {d^2+3d+2}6.$$

Plane irreducible curves of degree $d$ with $r$ nodes have genus
$$ g= \frac{d^2-3d+2}2-r = \frac{d^2}3-2d+\frac 23.
$$
But a general curve of genus $g$ has no linear series $g^2_d$. The Brill-Noether theory predicts that the moduli space of curves of genus $g$
contains a locus of dimension
$$(3g-3) -\rho(d,g,2)= (4g-3)-(3)(d-g-2)=\frac {d^2}3+d-\frac {25}3
$$
of curves with a series $g^2_d$, whose general element $\Gamma$ has a unique linear series $g^2_d$ which sends it to a plane curve of degree $d$ with $r$ nodes. 

If the guess is true, by \cite{Treger89} we get a family of dimension $\delta = (\frac {d^2}3+d-\frac {25}3)+8$ of sets $S$ of cardinality $r$ which are nodes of plane 
curves of degree $d$ (we added the automorphisms of $\p^2$). It is immediate to check that
$$ \dim (\p^2)^{(r)}-\delta = 2r-\delta = \frac {d^2+3d+2}3- \frac {d^2}3-d+\frac 13 =1. $$
\end{example}

The previous example does not apply when the degree $d$ is divisible by three, as Example \ref{noniche} will show.
\smallskip

\begin{remark}\label{tolgouno} 
Let $X$ be a variety of dimension $n$. We have for every $r$ a natural rational map $u_r:X^r\to X^{(r)}$. Moreover
we have a natural map $i:X^r\to X^{r-1}$ which forgets the first element.

Assume that $T$ is an irreducible  subvariety of codimension $1$ in $X^{(r)}$, not contained in the diagonal.
Then a general $S'\in X^{(r-1)}$ is contained in some $S\in T$. 
Moreover, the family $T_{S'}=\{S\in T: S'\subset S\}$ has dimension $n-1$. Indeed the set $u_r(i^{-1}(u^{-1}_{r-1}(S')))$, which parameterizes the sets
$S\in  X^{(r)}$ which contain $S'$, has dimension $n-1$ and intersects the diagonal in a finite set.
\end{remark}

\begin{example}\label{noniche} Let us see  what happens for the Veronese surface with $d=9$. Here the study of the Terracini locus makes sense only for $r\leq 18$.
Let us examine closely what happens for $r=18$.

Since $h^0(\Oo_{\p^2}(9))-3\cdot 18 =1$, by the AHThm for a general choice of a set of $18$ points $S$ there exists a unique nonic curve with $18$ 
nodes at the points of $S$.
On the other hand, there are special sets $S$ for which one finds a positive dimensional linear series of nonics passing through $2S$. Just to give an example,
if $S$ lies in a quartic curve $Q$, then for any line $R$ the nonic $2Q+R$ contains $2S$, so that $h^0(\Ii_{2S}(9))>0$ and 
$S$ lies in the Terracini locus $\t_{18}(\p^2,\Oo_{\p^2}(9))$.

We want to show that there are no components of $\t_{18}(\p^2,\Oo_{\p^2}(9))$ of codimension $1$ in $(\p^2)^{(18)}$.

Namely, by the AHThm and by \cite{CCi02a}, for a general set $S$  of $18$ points in $\p^2$ there exists a unique nonic curve which has nodes at $S$ and no other singularities,
and such a curve must be irreducible. 

If there exists a component $W$ of the Terracini locus of dimension $2\cdot 18-1= 35$, then 
the set of pairs  $\{(S,C): C$ is a 
nonic curve passing through $2S, S\in W\}$  produces a family of dimension $36$ of nodal nonics that cannot contain the general nonic of genus $10$.  
If the general member of the family is irreducible we get a contradiction with the irreducibility  of the family of nonics with $18$ nodes, proved in  \cite{Harr86}.

It remains to analyze what happens if the general nonic curve of the family is reducible. This can be done case by case.
Notice that, by Remark \ref{tolgouno}, if some component $T$ of $\t_{18}(\p^2,\Oo_{\p^2}(9))$ has codimension $1$, then a general subset $S'$ of $17$ points
in $\p^2$ sits in some element of $T$.

First assume that a general nonic plane curve $C$ containing $2S$, with $S\in T$, is non reduced. If $C$ is equal to a 
double line $2R$ plus a reduced curve $C'$ of degree $7$ then
$R$ cannot contain  more than $2$ points of $S'$. The remaining $15$ points of $S'$ must belong to $C'$. Thus $C'$ is singular at $15$ general points. Hence 
by the AHThm $C'$ moves in a family of dimension at most $H^0(\Oo_{\p^2}(7))-1-15=20$. Even adding all lines, we find that the family of curves $C$ 
has dimension at most $22$.
Among the nodes, only those in $R$ can move when we fix the curve $C$. We get then that the family of sets $S$ arising from curves like $C$ has dimension
at most $2+22<35$. With a similar trick one excludes all the cases in which $C$ is non reduced.

Assume that a general nonic plane curve $C$ containing $2S$, with $S\in T$, is reducible but reduced. If $C$ is equal to a line $R$ 
plus a curve $C'$ of degree $8$ then
$R$ cannot contain  more than $2$ points of $S'$. The remaining $15$ points of $S'$ must belong to $C'$. Thus $C'$ is singular at $15$ general points. 
Thus by the AHThm $C'$  moves
in a family of dimension at most $H^0(\Oo_{\p^2}(8)-1-15=29$. Even adding all lines, we find that the family of curves $C$ has dimension at most $31$.
No nodes can move when the curve $C$ is fixed. We get then that the family of sets $S$ arising from curves like $C$ has dimension
at most $2+29<35$. With a similar trick one excludes all the cases in which $C$ is reducible and  reduced.
\end{example}

Components of codimension $1$ of the Terracini locus are not easy to find, in general. We show that if the cardinality $r$
is small, then no such component exists in the case of Veronese varieties.

\begin{theorem}\label{cc3}
Fix integers $n\ge 2$, $d\ge 3$ and $r>0$ such that $r\le \binom{n+\lceil d/2\rceil -1}{n}-1$. Then $\dim \t_r(\p^n, \Oo_{\p^n}(d)) \le rn-n$.
In particular, no components of the Terracini locus have codimension $1$ in $(\p^n)^{(r)}$.
\end{theorem}
\begin{proof}
Take an irreducible component of the Terracini locus  $T$ with maximal dimension and assume $\dim T \geq rn-n+1$. By the AHThm
$T$ does not contain a general element of $(\p^n)^{(r)}$. 
Fix a general subset $S'\in(\p^n)^{(r-1)}$. By Remark \ref{tolgouno}, there exists a $(n-1)$-dimensional family of $S\in T$ such that $S'\subset S$. 
Thus the locus $E$ of points $o\in\p^2\setminus S'$ such that $S'\cup \{o\}\in T$ has dimension $n-1$. This implies that for a general line $R\subset\p^n$,
the intersection $R\cap E$ is non-empty (notice that $R\cap S'=\emptyset$ by the generality of $R$). We may assume $S=S'\cup\{o\}\in T$ for $o\in R$.

Let us prove that $h^1(\Ii_S(\lfloor d/2\rfloor))=0$. Since $S'$ is general and $\ell(S')< \binom{n+\lceil d/2\rceil -1}{n} $ we have $h^1(\Ii_{S'}(\lceil d/2\rceil -1))=0$.
Call $B$ the base locus of the linear system defined by $H^0(\Ii_{S'}(\lceil d/2\rceil -1))$. We claim that $\dim(B)\leq n-2$. Indeed by assumption
$\ell(S')\leq \binom{n+\lceil d/2\rceil -1}{n} -2$ so $h^0(\Ii_{S'}(\lceil d/2\rceil -1))\geq 2$. By the generality of $S'$, $S'$ sits in an integral hypersurface
$W$ in $H^0(\Ii_{S'}(\lceil d/2\rceil -1))$
and $B\neq W$. Since $B$ is closed, it follows $\dim B\leq n-2$. Thus $B\neq E$. Since $R$ is general,  $o$ is a general point of $E$, thus 
$o\notin B$. It follows that
$$ h^0(\Ii_{S}(\lceil d/2\rceil -1))<h^0(\Ii_{S'}(\lceil d/2\rceil -1))$$
which implies that $h^1(\Ii_{S'}(\lceil d/2\rceil -1))=0$. Hence $h^1(\Ii_S(\lfloor d/2\rfloor))=0$.

Next, we prove that $S$ satisfies $\dagger$ with respect to $\Oo_{\p^n}(\lceil d/2\rceil)$. We need to prove that for any $p\in S$ we have
$h^1(\Ii_{S\cup (2p,\p^n)}(\lceil d/2\rceil))=0$. 

Assume first $p=o$. Fix a general hyperplane $H$ containing $o$.
Since $H$ is general, we have $H\cap S'=\emptyset$ and hence 
$$H\cap (S\cup (2o,\p^n)) = H\cap (2o,\p^n) = (2o, H).$$ 
Moreover the residue of $S\cup (2o,\p^n)$ with respect to $H$ is $S$.  We get an exact sequence
\begin{equation}\label{eqcc2}
0\to \Ii _S(\lceil d/2\rceil -1)\to \Ii _{S\cup (2o,\p^n)}(\lceil d/2\rceil)\to \Ii _{(2o,H),H}(\lceil d/2\rceil)\to 0
\end{equation}
Since $\Oo_H(\lceil d/2\rceil)$ is very ample, we have $h^1(H,\Ii _{(2o,H),H}(\lceil d/2\rceil))=0$. Above we proved that $h^1(\Ii _S(\lceil d/2\rceil -1))=0$. Then
the long cohomology exact sequence of \eqref{eqcc2} gives $h^1(\Ii _{S\cup (2o,\p^n)}(\lceil d/2\rceil))=0$. 

Similarly, fix a point $p\in S'$ and let $M$ be a general
hyperplane containing $p$. As above $(S\cup (2p,\p^n))\cap M=(2p,M)$ and $Res(S\cup (2p,\p^n),M)=S$, so the exact sequence 
$$ 0\to \Ii _S(\lceil d/2\rceil -1)\to \Ii _{S\cup (2p,\p^n)}(\lceil d/2\rceil)\to \Ii _{(2p,M),M}(\lceil d/2\rceil)\to 0$$
proves, as above, that $h^1(\Ii _{S\cup (2p,\p^n)}(\lceil d/2\rceil))=0$. 

We conclude that $S$ satisfies $\dagger$ with respect to $\Oo_{\p^n}(\lceil d/2\rceil)$. 
By Theorem \ref{i1}, we are done.
\end{proof}

The previous result is not sharp. Indeed, we expect that finding examples of Terracini loci of codimension $1$, when $r$ is not close to the upper bound,
is hard in general.

For completeness, we state and prove the following proposition, where we introduce a criterion to guarantee that a set $S$ has the property that 
the addition of a general point does not create a set in the Terracini locus.
 
\begin{proposition}\label{a2}
Fix a reduced finite set $S\subset \p^n$, of cardinality $r$, which does not belong to $\t_r(\p^n,\Oo_{\p^n}(d-2))$. Then for all $p\in \p^n\setminus S$
we have that $S\cup p$ does not belong to the Terracini locus $\t_{r+1}(\p^n,\Oo_{\p^n}(d))$.
\end{proposition}
\begin{proof} Write $2S$ and $2p$ for $(2S,\p^n)$ and $(2p,\p^n)$ resp.
By the Castelnuovo-Mumford's regularity (\cite{Mumford}, lecture 14) the assumption  $h^1(\Ii _{2S}(d-2)) =0$ implies that 
homogeneous ideal of $2S$ is generated in degree $d-1$. 
Thus there is a degree $d-1$
hypersurface $T\subset \p^n$ containing $2S$ and with $p\notin T$. The residual exact sequence of $T$ gives the exact sequence
\begin{equation}\label{eqa2x}
0\to \Ii _{2p}(1)\to \Ii _{2S\cup 2p}(d) \to  \Ii _{(2S,T),T}(d)  \to 0
\end{equation}
Since $\Oo_{\p^n}(1)$ is very ample, we have $h^1(\Ii _{2p}(1))=0$. Since $h^1(\Ii _{2S}(d-2)) =0$, we also have $h^1(\Ii
_{2S}(d)) =0$. The restriction map $H^0(\Oo _{\p^n}(d)) \to H^0(\Oo_{2S})$, which is surjective, factors through $H^0(\Oo _{\p^n}(d)) \to H^0(T,\Oo _T(d))$,
so that $h^1(T,\Ii _{(2S,T),T}(d))  =0$.
Thus \eqref{eqa2x} gives $h^1(\Ii _{2S\cup 2p}(d)) =0$.
\end{proof}

\section{Terracini locus and identifiability}\label{ident}

Through this section, we assume $X\subset \p^N$, with $N\geq r(n+1)-1$, and $L=\Oo_X(1)$.

We show how the notion of Terracini locus is linked to the notion of {\it identifiability} of a point with respect to the variety $X$.

\begin{definition} For a point $u\in S_r(X)$ we say thet $u$ is $r$-identifiable if there exists a unique set $S\subset X$ of cardinality $r$ such that
$u$ belongs to the span of $S$ (see e.g. \cite{COtt12}).

We say that $u$ is algebraically $r$-identifiable if there are no positive dimensional families of linearly independent subsets $S\subset X$ of cardinality $r$ such that
$u$ belongs to the span of $S$. (see e.g. \cite{AmendolaRaneSturm18}).
\end{definition}

When $X$ is $r$-defective, then no sets $S\subset X$ of cardinality $r$ are algebraically identifiable. 

Identifiability and algebraic identifiability are strictly connected. Of course the former implies the latter, but there are cases in which they are equivalent 
(see e.g. \cite{BallC13}). When $X$ is a Segre or a Veronese variety, the test on the algebraic identifiability of a tensor $u\in S_r(X)$ can be performed 
 by using the following, well known fact (see \cite{AngeCVan18}, \cite{AngeCMazzon19} or \cite{Mazzon20}).

\begin{proposition}  Assume that $u$ is not algebraically $r$-identifiable for $X$, so that there exists a family $\Ss$ of linearly independent subsets $S$
of cardinality $r$ such that $u$ belongs to the span of $S$. Then every $S\in\Ss$ sits in the $r$-th Terracini locus of $X$.
\end{proposition}
\begin{proof} For any $S=\{p_1,\dots,p_r\}\in\Ss$ the abstract secant variety $AS_r(X)$ is smooth of dimension $r(n+1)-1$ at $(u,p_1,\dots,p_r)$
and it contains a positive dimensional subvariety $W$ which maps to $u$ in the projection $\pi$ to $\p^N$. All tangent vectors to $W$, which are also
tangent vectors to to $AS_r(X)$, are killed by the differential of $\pi$. Thus, the differential of $\pi$ has rank $<r(n+1)-1$. 
Since the image of the tangent  space to $AS_r(X)$ at $(u,p_1,\dots,p_r)$ contains
the span of the tangent spaces to $X$ at $p_1,\dots,p_r$, the claim follows.
\end{proof}

Thus, in order to prove that a set $p_1,\dots,p_r$ which spans $u$ is unique, it is necessary to prove that $S=\{p_1,\dots,p_r\}$ does
not belong to the Terracini locus. A way to do that, in the case of tensors, is illustrated in the algorithms of \cite{AngeCMazzon19} and \cite{Mazzon20}.
\smallskip

Finally, we give an example of elements in the Terracini locus of a Segre variety.

\begin{example} Let $X$ be the Segre embedding of $\p^3\times\p^3\times\p^3$ in $\p^{63}$. By Section 5 of \cite{COtt12}, $X$ is not $6$ defective, though
no points of $S_6(X)$ are identifiable. Fix $(u,p_1,\dots,p_6)\in AS_6(X)$. Considering $S=\{p_1,\dots,p_6\}$ as a subset of $\p^3\times\p^3\times\p^3$,
the natural projections $\pi_i:\p^3\times\p^3\times\p^3\to\p^3$ determine $3$ subsets $\pi_i(S)\in\p^3$.

Assume that the three sets $\pi_i(S)$ are projectively equivalent. By Remark 5.1 of \cite{COtt12}, $S$ belongs to an infinite
continuous family of elliptic normal curves in $X$. Every curve determines, by Proposition 5.2 of \cite{CCi06}, a second set $S'=\{q_1,\dots,q_r\}$
which spans $u$, and these sets clearly move in a positive dimensional family, because they are uniquely determined
by the choice of the elliptic curve. Such sets $S'$ form a positive dimensional subariety of the $6$-th Terracini locus of $X$.
\end{example}

With the notation of the previous example, we can launch

\begin{conjecture} The subsets $S\in X$ of cardinality $6$ such that $\pi_i(S)$ and $\pi_j(S)$ are projectively equivalent, for some fixed $i\neq j$, fill
a dense open set in the $6$-th Terracini locus of $X$.
\end{conjecture}

\end{document}